\documentclass[a4paper,12pt]{amsart}
\usepackage{amssymb,amsmath,hyperref,eufrak,bbm}

\author[Florent Benaych-Georges]{Florent Benaych-Georges}\address{Florent Benaych-Georges, LPMA,  UPMC Univ Paris 6, Case courier 188, 4, Place Jussieu, 75252 Paris Cedex 05, France, and CMAP, \'Ecole Polytechnique, route de Saclay, 91128 Palaiseau Cedex, France.} \email{florent.benaych@upmc.fr}
\title[Rectangular $R$-transform and rectangular spherical integrals]{Rectangular $R$-transform as the limit of rectangular spherical integrals}

\thanks{This work was partially supported by the \emph{Agence Nationale de la
Recherche} grant ANR-08-BLAN-0311-03.}
\date{\today}
\newcommand{\tta}{\theta}

\newcommand{\bbm}{\begin{bmatrix}}
\newcommand{\ebm}{\end{bmatrix}}
\newcommand{\bes}{\begin{equation*}}
\newcommand{\ees}{\end{equation*}}
\newcommand{\be}{\begin{equation}}
\newcommand{\ee}{\end{equation}}
\newcommand{\beqy}{\begin{eqnarray}}
\newcommand{\eeqy}{\end{eqnarray}}
\newcommand{\beq}{\begin{eqnarray*}}
\newcommand{\eeq}{\end{eqnarray*}}

\newcommand{\Var}{\operatorname{Var}}

\newcommand{\Diag}{\operatorname{diag}}

\newcommand{\Pro}{\mathbb{P}}

\newcommand{\Tr}{\operatorname{Tr}}

\newcommand{\ninf}{\underset{n\to\infty}{\longrightarrow}}
\newcommand{\nminf}{\underset{n,m\to\infty}{\longrightarrow}}

\newcommand{\one}{\mathbbm{1}}

\newcommand{\E}{\mathbb{E}}

\newcommand{\R}{\mathbb{R}}
\newcommand{\C}{\mathbb{C}}

\newcommand{\K}{\mathbb{K}}

\newcommand{\ud}{\mathrm{d}}

\newcommand{\pro}{probability }

\newcommand{\f}{\frac}
\newcommand{\ff}{\frac{1}}
\newcommand{\lf}{\left}
\newcommand{\ri}{\right}

\newcommand{\st}{such that }
\newcommand{\la}{\lambda}
\newcommand{\La}{\Lambda}

\newcommand{\ste}{\, ;\, }

\newcommand{\eps}{\varepsilon}
\newcommand{\arc}{\boxplus_{\la}}

\newcommand{\bxp}{\boxplus}

\newcommand{\bck}{\backslash}

\newtheorem{Th}{Theorem}[section]

\newtheorem{lem}[Th]{Lemma}

\newtheorem{rmq}[Th]{Remark}
\newtheorem{cor}[Th]{Corollary}

\newenvironment{pr}{\noindent {\it Proof. }}{\hfill$\square$}

\long\def\symbolfootnote[#1]#2{\begingroup
\def\thefootnote{\fnsymbol{footnote}}\footnote[#1]{#2}\endgroup}

\begin{document}
\maketitle

\begin{abstract}In this paper,   we connect   rectangular free probability theory and spherical integrals. We prove the analogue, for rectangular or square non-Hermitian matrices, of a result that Guionnet and Ma\"\i da proved for Hermitian matrices in \cite{alice-mylene05}. More specifically, we study the limit, as $n,m$ tend to infinity, of  
$\ff{n}\log\E\{\exp[\sqrt{nm}\theta X_n]\}$, where $\theta\in \R$, $X_n$ is the real part of  an entry of $U_n M_n V_m$,  $M_n$ is a certain $n\times m$   deterministic matrix and $U_n, V_m$ are independent Haar-distributed  orthogonal or unitary matrices with respective sizes $n\times n$, $m\times m$. We prove that when  the singular law of $M_n$ converges to a   \pro measure $\mu$,  for $\theta$ small enough, 
this limit actually exists and can be expressed with the rectangular $R$-transform of $\mu$. This gives an interpretation of this transform, which linearizes the rectangular free convolution, as the limit of a sequence of log-Laplace transforms. \end{abstract}


\section*{Introduction}In this article, we study the limit, as $n,m$ tend to infinity in such a way that $n/m$ tends to a limit $\la\in [0,1]$, of  
\bes  \ff{n}\log\E\{\exp[\sqrt{nm}\theta\,\Re(\Tr( E_nU_n M_n V_m))]\},\ees
where $\theta\in \R$, $M_n$ is a certain $n\times m$   deterministic matrix, $U_n, V_m$ are independent Haar-distributed  orthogonal or unitary matrices with respective sizes $n\times n$, $m\times m$, $E_n$ is an $m\times n$ elementary matrix (i.e. a matrix which entries are all zero, except one of them, which is equal to one) and $\Re(\cdot)$ denotes the real part.

The departure point of this study is the work of Collins, Zinn-Justin, Zuber, Guionnet, Ma\"\i da,   \'Sniady, Mingo and Speicher who proved, in the papers \cite{zz03, collinsIMRN03, alice-mylene05, cs06, cs07, mingo-piotr-collins-speicher07},  that under various hypotheses on some $n\times n$ matrices $A_n$ and $B_n$ and a positive exponent $\alpha$, the asymptotics  of  \bes\ff{n^\alpha}\log\E\{\exp[{n}\theta\Tr( B_nU_n A_n U_n^*)]\}\ees  are related to   free \pro theory.
 For example, it has been proved \cite[Th. 2]{alice-mylene05} that if the    spectral law ({i.e.} uniform distributions on eigenvalues)  of the self-adjoint  matrix $A_n$ converges to  a compactly supported \pro measure $\mu$ and $F_n=\Diag(1, 0, \ldots, 0)$, then for $\theta$ small enough, 
 \be\label{220809.09h53.2}\ff{n}\log\E\{\exp[{n}\theta\Tr( F_nU_n A_n U_n^*)]\}\ninf \f{\beta}{2}\int_0^{\f{2\theta}{\beta}} R_\mu(t)\ud t,\ee
 where $R_\mu$ is the so-called {\it $R$-transform} of $\mu$ and $\beta=1$ or $2$ according to wether we consider real or complex matrices. The $R$-transform is an integral transform of \pro measures on $\R$. Its main property is that it linearizes the additive free convolution $\bxp$, the binary operation on \pro measures on $\R$  which can be    defined by the fact that for $A$, $B$ large self-adjoint random matrices with   spectral laws tending to $\mu_A,\mu_B$  (as the dimension goes to infinity) and $U$ a Haar-distributed  orthogonal or unitary  matrix independent of $A$ and $B$,   the spectral law of $A+UBU^*$  tends to $\mu_A\bxp \mu_B$: the free convolution $\bxp$ can be thought as the analogue,  for the spectral laws of certain  large random matrices, of the classical convolution for real random variables. For all \pro measures $\mu, \nu$ on $\R$, we have \be\label{270809.midi}\qquad R_{\mu\bxp \nu}(t)=R_\mu(t)+R_\nu(t)\qquad\textrm{ (for $t$ in a neighborhood of zero)}.\ee Hence  in the analogy between  the free convolution and the classical one, the $R$-transform plays the role of the log-Laplace transform, and   \eqref{220809.09h53.2} gives a concrete sense to this: the $R$-transform (more specifically  its primitive, which also satisfies \eqref{270809.midi}), is the limit of  a certain sequence of  log-Laplace transforms.

 Let us now describe the content of our paper.  For each  $\la\in [0,1]$, another free convolution, denoted by $\bxp_\la$ and called the {\it rectangular free convolution with ratio $\la$},     does the same job as $\bxp$ for the {\it singular laws} ({i.e.} uniform distributions on singular values)   of   rectangular $n\times m$ random matrices which dimensions $n,m$ tend to infinity in such a way that $n/m$ tends to $\la$: for $n, m$ large integers \st $n/m\simeq \la$, for $A,B$ some $n\times m$ real or complex matrices  with singular laws  $\nu_A,\nu_B$  and $U,V$ Haar-distributed  orthogonal or unitary matrices independent of $A$ and $B$,   the singular  law of $A+UBV$ is  approximately $\nu_A\bxp_\la \nu_B$ (see  \cite{b09} or the introduction of \cite{b09R=C} for a more precise definition of $\bxp_\la$). Like the $R$-transform for $\bxp$ and the log-Laplace transform for the classical convolution, an integral transform linearizes $\bxp_\la$. It is called {\it the rectangular $R$-transform with ratio $\la$} and is denoted by $C^{(\la)}$: for all \pro measures $\mu, \nu$ on $[0,+\infty)$, we have \be\label{270809.miditrois}\qquad C^{(\la)}_{\mu\bxp_\la \nu}(t)=C^{(\la)}_\mu(t)+C^{(\la)}_\nu(t)\qquad\textrm{ (for $t$  in a neighborhood of zero)}.\ee The main result of the paper gives an interpretation of the rectangular $R$-transform (more specifically its primitive, which also satisfies \eqref{270809.miditrois}) as the limit of a sequence of log-Laplace transforms: we prove that if the singular law of $M_n$ tends to a \pro measure $\mu$ and $n/m$ tends to a limit $\la\in [0,1]$ as $n,m$ tend to infinity, then for $\theta$ small enough, for $E_n$ a sequence of $m_n\times n$ elementary matrices, 
 \be\label{220809.09h53.2.caroll}\ff{n}\log\E\{\exp[\sqrt{nm}\theta\, \Re(\Tr( E_nU_n M_n V_m))]\}\nminf\beta\int_0^{\f{\theta}{\beta}}\f{C_\mu^{(\la)} (t^2)}{t}\ud t.\ee

 Let us mention that free probability theory  has initially  been built  in the area of operator algebras and that concrete relations between {\it free} and {\it classical} probability theory, like the ones of \eqref{220809.09h53.2} and \eqref{220809.09h53.2.caroll}, are not that common. 

Let us also mention that expectations of the exponential of traces of  polynomials of constant matrices and uniform orthogonal random matrices, which    have been extensively studied in  physics and also other areas, like information theory,    are often called {\it spherical integrals}. 
See e.g.   \cite{zvonkin97, alice-ofer-2002, alice-saint-flour} and  the references above for the case of square matrices. Spherical integrals involving rectangular matrices are considered   in the papers   \cite{swrect.2003,ghaderipoor08,{kabashima08}}. The quantities studied in the paper \cite{kabashima08} of Kabashima are closely related to the spherical integral we study here: in Equation (8), Kabashima considers exactly the same spherical integral   as ours, with the same hypotheses (he supposes the singular law of $X$ to have a limit and $p/N$ to stay bounded), except that he supposes $\theta$ to be on the imaginary line,  whereas in our paper, $\tta$ is real. He is not giving an argument that qualifies as a mathematical proof, but  he gives an asymptotic formula  in Equation (9). With our vocabulary, this formula expresses  $$\lim_{n\to \infty}\ff{n}\log \E\{\exp i\Tr(E_nU_nM_nV_m)\}$$ as a saddle point value of a certain function. The actual computation of this saddle point is not easy, but his Equation (9)  has a structure which is quite close to our Equation \eqref{14.4.11}.

The paper is organized as follows. In Section \ref{220809.09h53}, we state  the main result of the paper, Theorem \ref{310709.11h}, and 
 discuss it. In Section \ref{prelim.R-tr.270809}, we recall the precise definition of the rectangular $R$-transform and prove a result of continuity of the map $(\la, \mu)\longmapsto C_\mu^{(\la)}$. At last,  Section \ref{5809.11h32} is devoted to the proof of Theorem \ref{310709.11h}, following  the ideas of the proof of \cite[Th. 2]{alice-mylene05}.

\section{Main result}\label{220809.09h53}

\subsection{Statement}
Let us consider, for all $n\geq 1$, an integer $m_n\geq n$ such that, as $n$ tends to infinity, $n/m_n$ tends to a limit $\la\in [0,1]$ and an $n\times m_n$ real or complex nonrandom matrix $M_n$ whose singular values are strictly bounded, uniformly in $n$, by a constant $K$ and such that, as $n$ tends to infinity, the singular law of $M_n$  converges weakly to a \pro measure that we shall denote by $\mu$. Let us define, for $\theta\in \R$, \be\label{200809.9h32}I_n(\theta)=\ff{n}\log\E\{\exp[\sqrt{nm_n}\theta\,\Re(\Tr( E_nU_n M_n V_n))]\},\ee
where $U_n, V_n$ are independent Haar-distributed  orthogonal or unitary (according to whether $M_n$ is real or complex)  matrices with respective sizes $n\times n$, $m_n\times m_n$ and $E_n$ denotes an $m_n\times n$  elementary matrix (i.e. a matrix which entries are all zero, except one of them, which is equal to one).

In the case where $\la=0$, we also suppose that there is $\alpha<2$ \st \be\label{200809.18h19}\textrm{for $n$ large enough,  }\quad m_n\leq n^\alpha.\qquad\qquad
\ee

\begin{rmq}\label{remark.notation.20.08.09}{\rm Let $\K$ be either $\R$ or $\C$ according to wether we consider real or complex matrices. $I_n(\theta)$ can also be considered as the Laplace transform of a certain scalar product estimated at a pair of independent random vectors, one of them being   a uniform random vector of the unit sphere of $\K^n$ and the other one being  the projection, on $\K^n$, of a uniform random vector of the unit sphere of $\K^{m_n}$. Indeed, let us  denote the singular values of $M_n$ by $\mu_{n,1}, \ldots, \mu_{n,n}$ and introduce (see \cite{H&J85}) some orthogonal or unitary matrices $P_n, Q_n$ with respective sizes $n\times n$, $m_n\times m_n$ \st
 $$M_n=P_n\bbm \mu_{n,1}&&&0&\cdots&0\\
&\ddots&&\vdots&&\vdots\\
&&\mu_{n,n}&0&\cdots&0\ebm Q_n.$$Let also, for each $n$, $(i_n, j_n)$ be the index of the non-null entry of $E_n$.  Then the $j_n$th  row (resp. $i_n$th column)  $u_n=(u_{n,1},\ldots, u_{n,n})$ (resp. $v_n=(v_{n,1}, \ldots, v_{n, m_n})^t$) of $U_nP_n$ (resp.  $Q_nV_n$) is uniformly distributed on the unit sphere of 
$\K^n$ (resp. $\K^{m_n}$) and 
 one has 
\be\label{5809.12h05}I_n(\theta)=\ff{n}\log\E\{\exp[\sqrt{nm_n}\theta\,\Re(\sum_{k=1}^n u_{n,k}\mu_{n,k}v_{n,k})]\}.\ee}\end{rmq}

The main result of the article  is the following one. 
\begin{Th}\label{310709.11h} Set $\beta=1$ or $\beta=2$ according to wether we consider real or complex matrices. The function $I_n$ converges uniformly on every compact subset of $(-\beta K^{-1},\beta K^{-1})$ to the function $$I(\theta)=\beta\int_0^{\f{\theta}{\beta}}\f{C_\mu^{(\la)} (t^2)}{t}\ud t,$$ where $C_\mu^{(\la)}$ denotes the rectangular $R$-transform of $\mu$ with ratio $\la$ (its definition is recalled in Section \ref{prelim.R-tr.270809} below).
\end{Th}
\begin{rmq}{\rm Note that  the function $C_\mu^{(\la)} $ is analytic on $(-K^{-2}, K^{-2})$  and vanishes at zero, so $I$ is actually well defined and analytic on $(-K^{-1}, K^{-1})$.}\end{rmq}

\subsection{Particular cases where the matrices $M_n$ are square ($\la=1$) or asymptotically flat ($\la=0$)}

Let us recall that the {\it $R$-transform}  of a \pro measure $\nu$ is the function $$R_\nu(z)=G_\nu^{-1}(z)-\ff{z}, \qquad\textrm{ for $G_\mu(z)=\int \f{\ud \mu(t)}{z-t}$}$$ 
 (the convention  we use is the one of the analytic approach to freeness \cite{hiai,agz09}, which is not exactly the one of the combinatorial approach   \cite{ns06}: $R_\nu^{\textrm{combinatorics}}(z)=zR_\nu^{\textrm{analysis}}(z)$). 
 
 Let $\mu_s$ be the symmetrization of $\mu$, defined by $\mu_s(A)=\f{\mu(A)+\mu(-A)}{2}$ for all Borel subset $A$ of $\R$, and $\mu^2$ be the push-forward of $\mu$ by the function $t\mapsto t^2$.
\begin{cor}
In the particular case where $\la=1$ (resp. $\la=0$), the limit $I$ of  $I_n$ can be expressed via the $R$-transform of $\mu_s$ (resp. $\mu^2$) in the following way $$I(\theta)=\beta\int_0^{\f{\theta}{2}}{R_{\mu_s} (t)}\ud t, \qquad\textrm{ (resp. $I(\theta)=\beta\int_0^{\f{\theta}{2}}t{R_{\mu^2} (t^2)}\ud t$.)}$$\end{cor}

\begin{pr} It suffices to prove that $C^{(1)}_\mu(t^2)=tR_{\mu_s}(t)$ and that $C^{(0)}_\mu(t)=tR_{\mu^2}(t)$. The second equation can be found in  \cite[Lem. 3.2 or Sect. 3.6]{b09}. The first equation  follows from the fact that for all $\la$, $C^{(\la)}_\mu=C^{(\la)}_{\mu_s}$ and from the fact that for all symmetric \pro measure $\nu$, by   \cite[Sect. 3.6]{b09}, $C^{(1)}_\nu(z^2)=zR_\nu(z)$.\end{pr}

\subsection{Possible extensions of Theorem \ref{310709.11h}}

\subsubsection{Cumulants point of view} For $\la\in [0,1]$, the {\it rectangular free cumulants with ratio $\la$} of $\mu$ have been defined in \cite[Sect. 3.4]{b09} (see also \cite[Sect. 2.2]{bg07c}): this is the sequence $(c_{2k}(\mu))_{k\geq 1}$ linked to the moments of $\mu$ by \cite[Eq. (4.1)]{fbg05.inf.div}. Recall also that for $X$ a bounded real random variable, the {\it classical cumulants}  of $X$ are the numbers $\operatorname{Cl}_k(X)$ defined by the formula $$\log \E (e^{zX})=\sum_{k\geq 1}\f{\operatorname{Cl}_k(X)}{k!}z^k.$$ Differentiating formally the convergence $I_n(\theta)\longrightarrow I(\theta)$, one would get the following ``classical cumulants interpretation" of the rectangular free cumulants with ratio $\la$: for all positive integers $k$, 
\be\label{29.9.10.1ch} c_{2k}(\mu)= \lim_{n\to\infty}\f{(nm_n)^k}{n}\f{\beta^{2k-1}}{(2k-1)!}\operatorname{Cl}_{2k}(\Re(\Tr( E_nU_n M_n V_n))).\ee This formula can be considered as a ``rectangular analogue" of \cite[Th. 4.7]{collinsIMRN03}. The author believes that \eqref{29.9.10.1ch} can be proved rigorously
with one of the following two methods.   One  of them would be to use the {\it Weingarten calculus},  developed by Collins and \'Sniady,  for the computation of the expectation of moments of the entries of the matrices $U_n$, $V_n$ (as in the proof of \cite[Th. 4.7]{collinsIMRN03}). The other one  would   rely on the extension of  Theorem \ref{310709.11h} to complex values of $\tta$ and notice that the functions in question there are analytic in $\theta$ (so that their convergence implies the one of their derivatives).    

\subsubsection{Case where $M$ is chosen at random} If $M_n$ is also chosen at random, independently of $U_n$ and $V_n$, and the expectation, in \eqref{200809.9h32}, is taken with respect to the randomness of $U_n,V_n$ and $M_n$, then Theorem \ref{310709.11h} stays true in certain cases (for example if $M_n$ is a standard  Gaussian matrix divided by $\sqrt{m_n}$), but it  can easily be seen that Theorem \ref{310709.11h} is not true in general anymore. However, if the expectation, in \eqref{200809.9h32}, only concerns the randomness of $U_n$ and $V_n$, then $I_n(\theta)$ is a random variable, and for certain sequences of random matrices $M_n$, more than Theorem \ref{310709.11h}   can be said about its convergence. An important example is given by the case where   $M_n=A_n+P_nB_nQ_n$ with $A_n, B_n$ deterministic matrices having limit singular laws $\mu_A, \mu_B$ and $P_n, Q_n$ Haar-distributed  orthogonal or unitary  matrices with respective sizes $n$ and $m_n$. In this case, one can prove (with technical additional hypotheses), that almost surely, 
\beq &\lim_{n\to\infty}\ff{n}\log\int e^{\sqrt{nm_n}\theta\, \Re(\Tr(E_nU_nM_nV_n))}\ud U_n\ud V_n=&\\   &\lim_{n\to\infty}\ff{n}\log\int e^{\sqrt{nm_n}\theta \,\Re(\Tr(E_nU_nA_nV_n))}\ud U_n\ud V_n+ \ff{n}\log\int e^{\sqrt{nm_n}\theta\, \Re(\Tr(E_nU_nB_nV_n))}\ud U_n\ud V_n.&\eeq 
To prove it, the ideas are the same as the ones of \cite[Sect. 6.1 and 6.2]{alice-mylene05}: first prove that the random variable  $$\ff{n}\log\int e^{\sqrt{nm_n}\theta \,\Re(\Tr(E_nU_nM_nV_n))}\ud U_n\ud V_n$$ concentrates around its mean (to do it, use \cite[Cor. 4.4.30]{agz09} instead of   \cite[Lem. 24]{alice-mylene05} in order to be allowed to consider the complex case) and then  the ideas of \cite[Sect. 6.2]{alice-mylene05}.
Since the singular law of $M_n$ converges almost surely to $\mu_A\bxp_\la \mu_B$, this gives a new proof of \eqref{270809.miditrois}. 

\subsubsection{Strong continuity property for the rectangular spherical integrals}  
One other way to extend this work, suggested to us by a referee, would be to use \eqref{octobre.2010.1} to prove a strong continuity property for the  rectangular spherical integrals.
in the spirit of Proposition 2.1 or Lemma 2.3 of \cite{mylene07}.

\section{Preliminaries about the rectangular $R$-transform}\label{prelim.R-tr.270809}

Let $\mu$ be a  \pro measure on the real line which support is contained in $[-K, K]$, with $K> 0$ (we do not suppose $\mu$ to be symmetric, how it was the case in the initial definition of the rectangular $R$-transform).
Let us define the generating function of the moments of $\mu^2$  $$M_{\mu^2}(z)=\int_{t\in\R}\f{t^2z}{1-t^2z}\ud \mu(t) = \int_{t\in\R}\ff{1-t^2z}\ud \mu(t)  -1\qquad(z\in [0,K^{-2})).$$
It can easily be proved that $M_{\mu^2}$ is nonnegative and non decreasing on $[0, K^{-2})$.
Let us define, for $\la\in [0,1]$, $T^{(\la)}(z)=(\la z+1)(z+1)$, and  
$$H^{(\la)}_\mu(z)= zT¨^{(\la)}( M_{\mu^2}(z)).$$ Then $H^{(\la)}_\mu$ defines an increasing analytic diffeomorphism {(in this paper, for $I$ an interval, we shall call an \emph{analytic function on $I$}  a function on $I$   which extends analytically to an open subset of $\C$ containing $I$)} 
from $ [0, K^{-2})$ onto the (possibly unbounded) interval $ [0, \lim_{z\uparrow K^{-2}}H^{(\la)}_\mu(z))$
 \st \bes\label{3809.11h38} H^{(\la)}_\mu(0)=0, \qquad
 \partial_z H^{(\la)}_\mu(0)=1, \qquad H^{(\la)}_\mu(z)\geq z, \qquad \lim_{z\uparrow K^{-2}}H^{(\la)}_\mu(z)\geq K^{-2}.\ees 
 We denote its inverse by ${H_\mu^{(\la)}}^{-1}$.
 Moreover, $T^{(\la)}$ defines an analytic increasing diffeomorphism from $[-1, +\infty)$ to $[0, +\infty)$, thus one can define  the {\em rectangular $R$-transform with ratio $\la$} of $\mu$:  
 \be\label{3809.11h39}C^{(\la)}_\mu(z)={T^{(\la)}}^{-1}\lf( \f{z}{{H_\mu^{(\la)}}^{-1}(z)}\ri)\textrm{ for $z\neq 0$,\quad and \quad$C^{(\la)}_\mu(0)=0$,}\ee which is analytic and non negative on 
the interval $ [0, \lim_{z\uparrow K^{-2}}H^{(\la)}_\mu(z))$ (which always contains $[0,K^{-2})$).  

By Theorems 3.8 and 3.12 of  \cite{b09}, the  rectangular $R$-transform characterizes symmetric measures, and for all pair $\mu_1, \mu_2$ of compactly supported symmetric \pro measures, $\mu_1\bxp_\la\mu_2$ is characterized by the fact that in a neighborhood of zero, \bes\label{28.08.08.17}
C^{(\la)}_{\mu_1\arc\mu_2}(z)=C^{(\la)}_{\mu_1}(z)+C^{(\la)}_{\mu_2}(z).\ees

The following theorem states the continuity of the mapping $(\la, \mu)\longmapsto C_\mu^{(\la)}$ in a way which is quite different from the one of Theorem 3.11 of  \cite{b09} (where $\la$ was fixed).
\begin{Th}\label{4809.19h09}Fix $K>0$, let $\mu_n$ be a sequence of  \pro measures on $[-K, K]$ which converges weakly to a limit $\mu$, and let $\la_n$ be a sequence of elements of $[0,1]$ which converges to a limit $\la\in [0,1]$. Then the sequence of functions $C^{(\la_n)}_{\mu_n}$ converges to $C^{(\la)}_{\mu}$ uniformly on every compact subset of $[0, K^{-2})$.
\end{Th}

\begin{pr} Recall that $C^{(\la)}_\mu$ is defined by \eqref{3809.11h39}. Since, by Heine's Theorem, $(\la,z)\mapsto {T^{(\la)}}^{-1} (z)$ is uniformly continuous on every compact subset of $[0,1]\times [0,+\infty)$, it suffices to prove that $\f{z}{{H_{\mu_n}^{(\la_n)}}^{-1}(z)}$ converges to $\f{z}{{H_{\mu}^{(\la)}}^{-1}(z)}$ uniformly on every compact subset of $[0, K^{-2})$. 
\\
\underline{Claim a} : {\it For each  compact subset $E$ of $\C\bck [K^{-2}, +\infty)$, there is a constant $k_E$ \st for any law $\nu$ on $[-K, K]$, for any $c\in [0,1]$,  for any $z\in E$, $$ |T^{(c)}(M_{\nu^2}(z))|\leq k_E.$$} Indeed, for  $z\in \C\bck [K^{-2}, +\infty)$, for any law $\nu$ on $[-K, K]$, for any $c\in [0,1]$, $$T^{(c)}(M_{\nu^2}(z))=\int_{(t,t')\in [-K,K]^2}\ff{(1-zt^2)(1-z{t'}^2)}\ud\nu(t)\ud(c\nu+(1-c)\delta_0)(t'),$$ thus $k_E=\max\{ |1-zt^2|^{-2}\ste |t|\leq K, z\in E\}$ is convenient. 
\\
\underline{Claim b} : {\it As $\nu$ varies in the set of laws on $[-K,K]$ and $c$ varies in $[0,1]$, the set of functions $$ z\in [0, K^{-2})\longmapsto\f{z}{{H_{\nu}^{(c)}}^{-1}(z)}  $$   is relatively compact for the topology of uniform convergence on every compact subset of $[0, K^{-2})$.}
   By Ascoli's Theorem, to prove Claim b, it suffices to prove that this family is uniformly bounded and uniformly Lipschitz on every compact subset of $[0, K^{-2})$. Let us fix $\nu$ a law on $[-K, K]$ and $c\in [0,1]$. Note that we have 
$$\f{z}{{H_{\nu}^{(c)}}^{-1}(z)}=\f{{H_{\nu}^{(c)}}(z)}{z}\circ{H_{\nu}^{(c)}}^{-1}(z), \quad\partial_z\f{z}{{H_{\nu}^{(c)}}^{-1}(z)}= \f{zH_\nu^{(c)'}(z)-H_\nu^{(c)}(z)}{z^2H_\nu^{(c)'}(z)}\circ{H_{\nu}^{(c)}}^{-1}(z).$$
Since, moreover, for all  $z\in [0, K^{-2})$, ${H^{(c)}_\nu}^{-1}(z)\leq z$ (indeed,  for all $z\in [0, K^{-2})$, $H^{(c)}_\nu(z)\geq z$), it suffices to verify that the sets of functions 
\beq &\{z\longmapsto \f{{H_{\nu}^{(c)}}(z)}{z}\ste \nu\textrm{ law on $[-K, K]$, }c\in [0,1]\}&\\ &\textrm{ and }\{ z\longmapsto\f{zH_\nu^{(c)'}(z)-H_\nu^{(c)}(z)}{z^2H_\nu^{(c)'}(z)}\ste \nu\textrm{ law on $[-K, K]$, }c\in [0,1]\}&\eeq   are uniformly bounded on every compact subset of $[0, K^{-2})$. The family of functions 
$\f{{H_{\nu}^{(c)}}(z)}{z}=T^{(c)}(M_{\nu^2}(z))$, indexed by $\nu,c$, is a family of  analytic functions on $\C\bck [K^{-2}, +\infty)$ which is uniformly bounded on every compact subset of $\C\bck [K^{-2}, +\infty)$ (by Claim a). As a consequence, the family of the derivatives $\partial_z \f{{H_{\nu}^{(c)}}(z)}{z}$ is  also uniformly bounded on every compact subset of $\C\bck [K^{-2}, +\infty)$. Since $$\f{zH_\nu^{(c)'}(z)-H_\nu^{(c)}(z)}{z^2H_\nu^{(c)'}(z)}=\ff{H_\nu^{(c)'}(z)}\partial_z \f{{H_{\nu}^{(c)}}(z)}{z}$$ and $H_\nu^{(c)'}(z)\geq 1$ on $[0, K^{-2})$, Claim b is proved.

Hence one can suppose that  $\f{z}{{H_{\mu_n}^{(\la_n)}}^{-1}(z)}$ converges to a function $f$ uniformly on every compact of $[0, K^{-2})$. Let us fix $z\in [0, K^{-2})$ and let us prove that $f(z)=\f{z}{{H_{\mu}^{(\la)}}^{-1}(z)}$. If $z=0$, it is clear (since all these functions are implicitly defined to map $0$ to $1$). Suppose that $z>0$. Note that $f(z)\neq 0$, because  for all $n$, $\f{z}{{H_{\mu_n}^{(\la_n)}}^{-1}(z)}\geq 1$.  Let us denote $l=\f{z}{f(z)}$. It suffices to prove that $l= {H_{\mu}^{(\la)}}^{-1}(z)$, i.e. that ${H_{\mu}^{(\la)}}(l)=z$. Since $${H_{\mu}^{(\la)}}(l)=\lim_{n\to\infty}{H_{\mu}^{(\la)}}({H_{\mu_n}^{(\la_n)}}^{-1}(z)),$$ it suffices to prove that  $H_{\mu_n}^{(\la_n)}$  converges to $H_{\mu}^{(\la)}$ uniformly on every compact subset of $[0, K^{-2})$. 
But it is  easy to see, using \cite[Th. C.11]{agz09}, that $M_{\mu_n^2}$ converges to $M_{\mu^2}$ uniformly on every compact subset of $[0, K^{-2})$ and then that $H_\mu^{(\la_n)}$  converges to $H_{\mu}^{(\la)}$ uniformly on every compact subset of $[0, K^{-2})$. 
\end{pr}

\section{Proof of Theorem \ref{310709.11h}}\label{5809.11h32}
\subsection{Preliminaries}
We shall use the following lemmas  several times in  the paper. Let $\|\cdot\|$ denote the canonical euclidian norm on each $\R^d$. 
\begin{lem}\label{31709.17h56}
Let $(G_i)_{i\geq 1}$ be a family of independent real random variables with standard Gaussian law. Let $T$ be fixed  and let, for each $n$, $(\sigma_{n,1}, \ldots, \sigma_{n,n})\in [0,T]^n$ be \st  \be\label{31709.17h57}\ff{n}\sum_{i=1}^n\sigma_{n,i}^2=1.\ee Let us define, for each $n$,  $X_n=(\sigma_{n,1}G_1, \ldots, \sigma_{n,n}G_n)$. Then for all $\kappa\in (0, \ff{2})$, for $n$ large enough, \bes\label{31709.17h55}\Pro\{|\|X_n\|-\sqrt{n}|\geq n^{\ff{2}-\kappa} \}\le 2T^4n^{2\kappa-1}.\ees 
\end{lem}

\begin{pr}
Note that by \eqref{31709.17h57}, the random variable $N_n:=\f{\|X_n\|^2}{n}-1$ is centered. Moreover, $\Var(N_n)=\f{\Var (G_1^2) }{n^2}\sum_{i=1}^n\sigma_{n,i}^4\leq \f{2T^4  }{n}$. It follows, by  Tchebichev's inequality, that for all  $\kappa\in (0, \ff{2})$, $$\Pro \{|N_n|\geq n^{-\kappa}\}\leq 2T^4 n^{2\kappa-1}.  $$ To deduce that for $n$ large enough,$$\Pro \lf\{\lf|\f{\|X_n\|}{\sqrt{n}}-1\ri|\geq n^{-\kappa}\ri\}  \leq 2T^4 n^{2\kappa-1},  $$ it suffices to notice that the function $\sqrt{\cdot}$ is 1-Lipschitz on $[1/4, +\infty)$ and that $n^{-\kappa}\leq 3/4$ for $n$ large enough. 
\end{pr}

\begin{lem}
Let $\mu$ be a  \pro measure whose support is contained in $[-K, K]$, fix $\la\in [0,1]
$, $\theta\in [0, K^{-1})$ and define $\gamma= C_\mu^{(\la)}(\theta^2)$. Then \be\label{38.09.19h}M_{\mu^2}\lf(\f{\theta^2}{T^{(\la)}(\gamma)}\ri)=\gamma.\ee 
\end{lem}

\begin{pr}By the definition of $C^{(\la)}_\mu$ given in \eqref{3809.11h39}, $\f{\theta^2}{{H_\mu^{(\la)}}^{-1}(\theta^2)}=T^{(\la)}(\gamma)$, hence $\f{\theta^2}{T^{(\la)}(\gamma)}={H_\mu^{(\la)}}^{-1}(\theta^2)$. Since $\gamma\geq 0$, $\f{\theta^2}{T^{(\la)}(\gamma)}\in [0, K^{-2})$ and one can apply the function ${H_\mu^{(\la)}}$ on both sides. We get ${H_\mu^{(\la)}}\lf(\f{\theta^2}{T^{(\la)}(\gamma)}\ri)=\theta^2$, i.e. $$\f{\theta^2}{T^{(\la)}(\gamma)}T^{(\la)}\lf( M_{\mu^2}\lf(\f{\theta^2}{T^{(\la)}(\gamma)}\ri)\ri)=\theta^2 .$$ It follows that $ T^{(\la)}\lf(M_{\mu^2}\lf(\f{\theta^2}{T^{(\la)}(\gamma)}\ri)\ri)={T^{(\la)}(\gamma)}$. Since both $M_{\mu^2}\lf(\f{\theta^2}{T^{(\la)}(\gamma)}\ri)$ and $\gamma$ are nonnegative real numbers, one gets 
\eqref{38.09.19h}.
\end{pr}

\begin{lem}\label{lastcarol.momo}Let $X_n$ be a sequence of nonnegative random variables, with positive expectations. Let $Z_n$ be a sequence of real random variables \st there exists deterministic constants $C, \eta>0$ \st  for all $n$, $|Z_n|\leq Cn^{1-\eta}$. Then as $n$ tends to infinity, $$\ff{n}\log \E(X_ne^{Z_n})=\ff{n}\log \E(X_n)+o(1).$$
\end{lem}

\begin{pr}It suffices to notice that we have $X_ne^{-Cn^{1-\eta}}\leq X_ne^{Z_n}\leq X_ne^{Cn^{1-\eta}}$.
\end{pr}

\subsection{Notation for the proof of Theorem \ref{310709.11h}} In the next sections, $o(1)$ shall denote any sequence  of functions on $(-K^{-1}, K^{-1})$ which converges to zero as $n$ tends to infinity, uniformly on every compact subset of $(-K^{-1}, K^{-1})$. Also, we shall work with the notation introduced in Remark \ref{remark.notation.20.08.09} and handle $I_n(\theta)$ via Formula \eqref{5809.12h05}. 

\subsection{Reduction to the case where all singular values of $M_n$ are positive}Let us suppose the result to be proved in the particular case where for all $n,k$, $\mu_{n,k}>0$, and let us prove it in the general case. We set, for each $n,k$,  $$\tilde{\mu}_{n,k}=\mu_{n,k}+\min \{ m_n^{-2}, ({K-\mu_{n,k}})/{2}\}$$
and define the perturbation of $I_n(\tta)$ 
\bes\label{5809.12h05_sufjan}\tilde{I}_n(\theta)=\ff{n}\log\E\{\exp[\sqrt{nm_n}\theta\sum_{k=1}^n \Re(u_{n,k}\tilde{\mu}_{n,k}v_{n,k})]\}.\ees The uniform law on the $\tilde{\mu}_{n,k}$'s converges weakly to $\mu$ as $n$ tends to infinity and we have $0<\tilde{\mu}_{n,k}<K$, so by hypothesis, if follows that 
  $\lim_{n\to\infty} \tilde{I}_n(\theta)=I(\tta)$. Note moreover that   $$\tilde{I}_n(\theta)=\ff{n}\log\E\{\exp[\sqrt{nm_n}\theta\sum_{k=1}^n\Re( u_{n,k}{\mu}_{n,k}v_{n,k})]e^{Z_n}\},$$ with $$Z_n= \sqrt{nm_n}\theta \sum_{k=1}^n \Re(u_{n,k}v_{n,k})\min \{ m_n^{-2}, ({K-\mu_{n,k}})/{2}\}.$$
Since $|Z_n|\le K^{-1}$,      Lemma \ref{lastcarol.momo} allows us to claim that  $\lim_{n\to\infty}  {I}_n(\theta)=I(\tta)$.

\subsection{Deducing the complex case from the real one} Let   us explain how one can deduce the complex case from the real one. Let us use, in this paragraph, the notation $I^{(\beta)}_{n,m_n}(\tta,M_n)$ to emphasize on the value of each of the parameters defining $I_n(\tta)$. We have \beqy I^{(2)}_{n,m_n}(\tta,M_n)&=&\ff{n}\log\E\{\exp[\sqrt{nm_n}\theta\sum_{k=1}^n \Re(u_{n,k}{\mu}_{n,k}v_{n,k})]\} \nonumber \\
 &=&\ff{n}\log\E\{\exp[\sqrt{nm_n}\theta\sum_{k=1}^n {\mu}_{n,k}(\Re(u_{n,k})\Re(v_{n,k})-\Im(u_{n,k})\Im(v_{n,k}))]\} \nonumber \\
&=&2I_{2n,2m_n}^{(1)}\lf(\f{\tta}{2}, \bbm M_n& 0\\0& -M_n\ebm\ri)\label{22.9.10J-1}
\eeqy
Indeed, for $(z_{n,1}, \ldots, z_{n,n})$ a vector with uniform distribution on the unit sphere of $\C^n$ the vector $(\Re(z_1), \Im(z_1), \ldots, \Re(z_n), \Im(z_n))$ has uniform distribution on the unit sphere of    $\R^{2n}$ (use the realization of such vectors as the projection of standard Gaussian vectors on the sphere to see it). 

Since the singular values of $\bbm M_n& 0\\0& -M_n\ebm$ are the one of $M_n$ with multiplicity multiplied by two,   \eqref{22.9.10J-1} allows to deduce the complex case from the real one.

Let us now prove Theorem \ref{310709.11h} in the case where all $\mu_{n,k}$'s are positive and where $\beta=1$.  
\subsection{Proof of Theorem \ref{310709.11h}: a) Expression of $I_n(\theta)$ with a Gaussian integral} As written  above, we suppose from now on that all $\mu_{n,k}$'s are positive and that $\beta=1$.  In this section, we shall first explain how to replace, in Formula \eqref{5809.12h05}, $\sqrt{n} u_{n,k}$
and $\sqrt{m_n}v_{n,k}$ by  independent standard Gaussian 
variables (Formula \eqref{210809.10h56prime}) and then inject (quite artificially first) $C_\mu^{(\la)}(\theta^2)$ in the formula of $I_n(\theta)$ (Equation \eqref{bird.31709}).

For each $n$, let us define the function $$f_n : ((x_1, \ldots, x_n),(y_1, \ldots, y_{m_n}))\in \R^n\times \R^{m_n}\longmapsto \sum_{k=1}^n \mu_{n,k}x_ky_k.$$

By \cite[Sect. 1.2.1]{alice-mylene05}, up to a change of the \pro space which does not change the expectation, one can suppose that there are independent standard Gaussian vectors $x_n,y_n$ of  $\R^n$, respectively $ \R^{m_n}$,  
 \st $$u_n=\f{x_n}{\|x_n\|}, \qquad v_n=\f{y_n}{\|y_n\|}.$$

Let us fix $\kappa\in (0,1/2)$. If $\la>0$, the precise choice of $\kappa\in (0,1/2)$ is irrelevant, but   
if $\la=0$,  we choose $\kappa\in( \f{\alpha-1}{2},\ff{2})$ ($\alpha$ is the one of \eqref{200809.18h19}). Let us now  define the set $$A_n:=\lf\{(x,y)\in \R^n\times \R^{m_n}\ste \lf|{\|x\|}-{\sqrt{n}}\ri|\leq n^{\ff{2}-\kappa},\lf|{\|y\|}-{\sqrt{m_n}}\ri|\leq m_n^{\ff{2}-\kappa}\ri\}.$$
The event 
$\{(x_n, y_n)\in A_n\}$ is   independent of $(u_n,v_n)$ (because the density of the law of a standard Gaussian vector is a radial function), thus 
 \beq I_n(\theta) &=&\ff{n}\log \E\lf[\one_{(x_n,y_n)\in A_n}\exp (\sqrt{nm_n}\theta f_n(u_n,v_n))\ri]-\ff{n}\log \Pro(A_n).\eeq Moreover,  by Lemma \ref{31709.17h56}, $\Pro\{(x_n, y_n)\in A_n\}\longrightarrow 1$ as $n\to\infty$, thus \be\label{miracle.0809.10h57} I_n (\theta) =\ff{n}\log \E\lf[\one_{(x_n,y_n)\in A_n}\exp (\sqrt{nm_n}\theta f_n(u_n,v_n))\ri]+o(1).\ee

Moreover, note that on the event 
$\{(x_n, y_n)\in A_n\}$, \beq \sqrt{n}-n^{\ff{2}-\kappa}\leq &\| x_n\|&\leq \sqrt{n}+n^{\ff{2}-\kappa}\\
\sqrt{m_n}-m_n^{\ff{2}-\kappa}\leq &\| y_n\|&\leq \sqrt{m_n}+m_n^{\ff{2}-\kappa},\eeq thus, since $m_n\geq n$,  \be\label{200809.21h42}
\sqrt{nm_n}-3\sqrt{m_n}n^{\ff{2}-\kappa} \leq \| x_n\| \| y_n\|\leq \sqrt{nm_n}+3\sqrt{m_n}n^{\ff{2}-\kappa}.\ee If $\la >0$, since $m_n/n$ is bounded, it follows that 
 there is a deterministic constant $C$ independent of $n$ \st on  on the event 
$\{(x_n, y_n)\in A_n\}$, \be\label{210809.10h55} |\,\|x_n\|  \|y_n\|-\sqrt{nm_n}\,|\leq C n^{1-\kappa}.\ee 
If $\la=0$, it follows from \eqref{200809.21h42} and \eqref{200809.18h19} that for $\eta=\f{1-\alpha}{2}+\kappa$ (which is positive by definition of $\kappa$),   for $n$ large enough, \be\label{210809.10h56} |\,\|x_n\|  \|y_n\|-\sqrt{nm_n}\,|\leq 3 n^{1-\eta}.\ee

 Note that by \eqref{miracle.0809.10h57}, \bes I_n(\theta)=\ff{n}\log \E[\one_{(x_n,y_n)\in A_n}e^{\theta f_n(x_n,y_n)+\f{\theta  f_n(x_n,y_n)}{\|x_n\| \|y_n\|}(\sqrt{nm_n}-\|x_n\| \|y_n\|)
} ]+o(1),\ees and that for all $n,k$,  $0\le \mu_{k,n}\leq K$, which implies that $\lf|\f{f_n(x_n,y_n)}{\|x_n\| \|y_n\|}\ri|\leq K$. Hence by Lemma \ref{lastcarol.momo} and \eqref{210809.10h55} (or \eqref{210809.10h56} if $\la=0$), \be\label{210809.10h56prime}I_n(\theta)=\ff{n}\log \E\lf[\one_{(x_n,y_n)\in A_n}e^{\theta f_n(x_n,y_n)} \ri]+o(1).\ee

Note that  on the event 
$\{(x_n, y_n)\in A_n\}$, we have
 \beq  n-2n^{1-\kappa}\leq n-2n^{1-\kappa}+n^{1-2\kappa}\leq &\| x_n\|^2&\leq n+2n^{1-\kappa}+n^{1-2\kappa}\leq n+3n^{1-\kappa}\\  
 n-2n^{1-\kappa}\leq
n-2nm_n^{-\kappa}+nm_n^{-2\kappa}\leq &\f{n}{m_n}\| y_n\|^2&\leq n+2nm_n^{-\kappa}+nm_n^{-2\kappa}\leq n+3n^{1-\kappa}.
\eeq 
thus  
  for all $n$,  on the event 
$\{(x_n, y_n)\in A_n\}$, \be\label{210809.11h20}\lf|\, \|x_n\|^2-n\,\ri|+  \big|\,\f{n}{m_n} \|y_n\|^2-n\,\big| \leq 6 n^{1-\kappa}.\ee
 Now, let us define, for each $n$, 
\be\label{5809.7h48}\gamma_n (\theta)= C_{\mu_n}^{(\f{n}{m_n})}(\theta^2)\qquad\textrm{ for $\mu_n=\ff{n}\sum_{k=1}^n \delta_{\mu_{n,k}}$.}\ee
 Note that $\mu_n$ is the singular law of $M_n$, which tends to $\mu$. Hence by Theorem \ref{4809.19h09},  we have
\be\label{5809.7h32} \gamma_n (\theta)\ninf  C_{\mu}^{(\la)}(\theta^2)\qquad\textrm{ uniformly on every compact subset of $(-K^{-1}, K^{-1})$,}\ee
 so by \eqref{210809.11h20}, for every such compact set $E$, 
there is a constant $Q_E$ \st for all $n$, for all $\theta\in E$,  on the event 
$\{(x_n, y_n)\in A_n\}$, we have $$  |\gamma_n (\theta)( \ff{2} \|x_n\|^2+\f{n}{2m_n}\|y_n\|^2-n)|  \leq Q_En^{1-\kappa}.$$ Hence, by \eqref{210809.10h56prime} and   Lemma \ref{lastcarol.momo}, \beq I_n(\theta)\!\!\!\!\!&=&\!\!\!\!\!\ff{n}\log \E\lf[\one_{(x_n,y_n)\in A_n}\exp\lf\{\theta f_n(x_n,y_n)-\gamma_n (\theta)\lf( \ff{2} \|x_n\|^2+\f{n}{2m_n}\|y_n\|^2-n\ri)\ri\} \ri]+o(1)\\ 
\!\!\!\!\!&=&\!\!\!\!\! \gamma_n (\theta)+ \ff{n}\log \underbrace{
\E\lf[
\one_{(x_n,y_n)\in A_n}\exp
\lf\{\theta f_n(x_n,y_n)-\gamma_n (\theta)\lf( \ff{2} \|x_n\|^2+\f{n}{2m_n}\|y_n\|^2\ri)
\ri\}
\ri]
}_{:=J_n(\theta)}+o(1).\eeq

Thus, by  \eqref{5809.7h32},    \be\label{bird.31709}I_n(\theta)= C_\mu^{(\la)}(\theta^2)+ \ff{n}\log J_n(\theta)+o(1).\ee

\subsection{Proof of Theorem \ref{310709.11h}: b) Convergence of the Gaussian integral}
We have, assimilating the vectors of $\R^n$ and $\R^{m_n}$ with column-matrices,  
\be\label{5809.19h}J_n(\theta)=(2\pi)^{-\f{n+m_n}{2}}\int_{x\in\R^n, y\in \R^{m_n}}\one_{A_n}(x,y)\exp\lf\{ -\ff{2}\bbm x^t&y^t\ebm T_n\bbm x\\ y\ebm \ri\}\ud x\ud y,\ee
 for $$T_n:=\begin{bmatrix}a_n(\theta)I_n& \La_n(\theta) &0_{n,m_n-n}\\ \La_n(\theta)&b_n(\theta)I_n&0_{n,m_n-n} \\ 0_{m_n-n,n}& 0_{m_n-n,n}&  b_n(\theta)I_{m_n-n} \end{bmatrix},$$ where $a_n(\theta)=1+\gamma_n(\theta)$, $b_n(\theta)=1+\f{n}{m_n}\gamma_n(\theta)$ and $\La_n(\theta)$ is the diagonal $n\times n$ matrix with diagonal entries $$\la_{n,1}(\theta):=-\theta\mu_{n,1},\ldots, \la_{n,n}(\theta):=-\theta\mu_{n,n}.$$ 

\noindent{\bf Notation:} In this section, in order to lighten the notation, we shall write $J_n$ for $J_n(\theta)$, $a_n$ for $a_n(\theta)$, {\it etc.} We shall also use the notation of the matricial  functional calculus   (thus assimilate $a$ and $aI_n$, {etc.}).
\begin{lem}\label{5809.19h46}
Let us fix $n\geq 1$ and let $a,b$ be real numbers and $\La$ an invertible  diagonal real $n\times n$ matrix.  Let us define    $$\Delta=(b-a)^2+4\La^2, \quad r^\pm=\f{a+b\pm \sqrt{\Delta}}{2}, \quad f^\pm =\ff{\sqrt{2\Delta\pm 2(b-a)\sqrt{\Delta}}}$$ and $$T=\begin{bmatrix}a& \La\\ \La& b \end{bmatrix}, \;\;D=\begin{bmatrix}r^+& 0\\ 0& r^-\end{bmatrix}, \;\;P=\begin{bmatrix}2\La f^+& 2\La f^-\\ (b-a)f^++\sqrt{\Delta}f^+& (b-a)f^--\sqrt{\Delta}f^- \end{bmatrix}. $$ Then $P$ is an orthogonal matrix and we have $T=PDP^t .$
\end{lem}

\begin{pr} One can easily verify that $P$ is orthogonal. Let us define $$Q=\begin{bmatrix}2\La & 2\La \\ (b-a)+\sqrt{\Delta}& (b-a)-\sqrt{\Delta} \end{bmatrix},H=\begin{bmatrix}f^+ & 0 \\ 0& f^- \end{bmatrix}.$$Then $P=QH$. One can easily verify that $TQ=QD$. If follows that $TQH=QDH$.  Since $HD=DH$,  $ TQH =QHD$, i.e. $TP=PD$, thus $T=PDP^t$.
\end{pr}

For $\theta\neq 0$, let us define $\Delta_n$, $r^\pm_n  $, $f_n^\pm $ as in the lemma, using $\La_n$ instead of $\La$, $a_n$ instead of $a$ and $b_n$ instead of $b$.  Let us define  $P_n$ in the same way, extended to an $(n+m_n)\times (n+m_n)$ matrix by adding $I_{m_n-n}$ on the  lower-right corner, i.e. $$P_n= \begin{bmatrix}2\La_n f_n^+& 2\La_n f_n^-&0\\ (b_n-a_n)f_n^++\sqrt{\Delta_n}f_n^+& (b_n-a_n)f_n^--\sqrt{\Delta_n}f^-&0\\ 0&0&I_{m_n-n} \end{bmatrix},$$and $D_n$ extended to an  $(n+m_n)\times (n+m_n)$ matrix by adding $b_nI_{m_n-n}$ on the lower-right corner, i.e. $$D_n=\begin{bmatrix}r_n^+& 0&0\\ 0& r_n^-&0\\ 0&0&b_nI_{m_n-n}\end{bmatrix} .$$ 
For $\theta=0$, we set $r_n^\pm=1$, $P_n=D_n=I_{n+m_n}$. 

Let us denote, for $X$ an $(n+m_n)\times (n+m_n)$ matrix,  $X(A_n)=\{X\bbm x \\ y\ebm \ste (x,  y) \in A_n\}$.  Let us also introduce  a standard Gaussian random column vector in $\R^{n+m_n}$, that we shall denote by  $$Z_n=(\underbrace{Z^+_{n,1},\ldots, Z^+_{n,n}}_{:=Z_n^+},\underbrace{Z^-_{n,1},\ldots, Z^-_{n,n}}_{:=Z_n^-},\underbrace{Z_{n,1}^0,\ldots,Z_{n,m_n-n}^0}_{:=Z_n^0})^t.$$ 
We have, by \eqref{5809.19h} and Lemma \ref{5809.19h46}, \beqy \nonumber
 J_n
 &=&
 (2\pi)^{-\f{n+m_n}{2}}\int_{A_n}\exp\{ -\ff{2}\bbm x^t & y^t\ebm P_nD_nP^t_n\bbm x\\ y\end{bmatrix}\}\ud x\ud y.\eeqy Thus, since $P_n$ is an orthogonal matrix, \beqy \nonumber
 J_n
 &=& (2\pi)^{-\f{n+m_n}{2}}\int_{P_n^t(A_n)}\exp\{ -\ff{2}\bbm x^t & y^t\ebm D_n\bbm x\\ y\end{bmatrix}\}\ud x\ud y.\eeqy Hence, by definition of $D_n$, we have 
\beqy\nonumber J_n&=& (2\pi)^{-\f{n+m_n}{2}}[b_n^{m_n-n}\prod_{i=1}^nr_{n,i}^+r_{n,i}^-]^{-1/2} \int_{\sqrt{D_n}P_n^t(A_n)}\!\!\!\!\!\!\!\!\!\!\!\exp\{ -\ff{2}(\| x\|^2+ \|y\|^2)\}\ud x\ud y,\eeqy which,  by definition of $Z_n$, can be written 
 \be\label{5809.8h15}
 J_n=  [b_n^{m_n-n}\prod_{i=1}^n(a_nb_n-\la_{n,i}^2)]^{-1/2}\underbrace{\Pro \{Z_n\in \sqrt{D_n}P^t_n(A_n)\}}_{=\Pro \{P_nD_n^{-1/2}Z_n\in A_n\}}
\ee

Let
 $X_n=(X_{n,1},\ldots,X_{n,n})^t$ be the vector of the first $n$ coordinates of $P_nD_n^{-1/2}Z_n$ and $Y_n=(Y_{n,1},\ldots,Y_{n,m_n})^t$
  be the one of the $m_n$ last ones. By definition of the set  $A_n$, we have \be\label{29.9.10} P_nD_n^{-1/2}Z_n\in A_n\iff \lf|{\|X_n\|}-{\sqrt{n}}\ri|\leq n^{-\kappa}\textrm { and }\lf|{\|Y_n\|}-{\sqrt{m_n}}\ri|\leq m_n^{-\kappa}.\ee 
  
  \underline{Claim} : {\it Both events 
$\lf\{\lf|{\|X_n\|}-{\sqrt{n}}\ri|\leq n^{-\kappa}\ri\}$ and $\lf\{\lf|{\| Y_n\|}-{\sqrt{m_n}}\ri|\leq m_n^{-\kappa}\ri\}$ have probabilities tending to one as $n$ tends to infinity,  
uniformly on every compact subset of $(-K^{-1}, K^{-1})$ (the random vectors $X_n$ and $Y_n$ depend indeed   on $\theta$).}

 For $\theta=0$, $X_n$ (resp. $Y_n$) is a standard Gaussian random vector of $\R^n$ (resp. $\R^{m_n}$).  For $\theta\neq 0$, by the definitions of $P_n$ and $ D_n$, 
 \beq X_n&=&2\La_n(f_n^+(r_n^+)^{-1/2}Z_n^++f_n^-(r_n^-)^{-1/2}Z_n^-),\\
Y_n&=&\begin{bmatrix}(b_n-a_n+\sqrt{\Delta_n})f_n^+(r_n^+)^{-1/2}Z_n^++(b_n-a_n-\sqrt{\Delta_n})f_n^-(r_n^-)^{-1/2}Z_n^-\\ b_n^{-1/2}Z_{n}^0\end{bmatrix}.
\eeq
Thus for each $n$, $X_n$ (resp. $Y_n$) has the law of $$(\sigma_{n,1}G_1, \ldots,\sigma_{n,n}G_n)\qquad\textrm{ (resp. }(\sigma'_{n,1}G_1, \ldots,\sigma'_{n,m_n}G_{m_n})\textrm{),}$$  for $(G_i)_{i\geq 1}$ a family of independent real random variables with standard Gaussian law and where if $\theta=0$, all $\sigma_{n,i}$'s and $\sigma_{n,i}'$'s are equal to $1$, whereas if   $\theta\neq 0$, for each $i=1,\ldots,n$,  \beq \sigma_{n,i}^2&=& 4\la_{n,i}^2\f{2}{ (2\Delta_{n,i}+2(b_n-a_n)\sqrt{\Delta_{n,i}} )( a_n+b_n+\sqrt{\Delta_{n,i}}) }\\ &&\qquad\qquad\qquad+ 4\la_{n,i}^2\f{2}{ (2\Delta_{n,i}-2(b_n-a_n)\sqrt{\Delta_{n,i}} )( a_n+b_n-\sqrt{\Delta_{n,i}})},\eeq 
 for each $i=n+1,\ldots, m_n$,
${\sigma'_{n,i}}^2=
b_n^{-1}$ and for each $i=1,\ldots,n$, \bes {\sigma'_{n,i}}^2=\f{2(b_n-a_n+\sqrt{\Delta_{n}})^2}{(a_n+b_n+\sqrt{\Delta_n})[2(b_n-a_n)\sqrt{\Delta_n}+2\Delta_n]} +\f{2(\sqrt{\Delta_{n}}-(b_n-a_n))^2}{(a_n+b_n-\sqrt{\Delta_n})[-2(b_n-a_n)\sqrt{\Delta_n}+2\Delta_n]}.\ees 

Hence by Lemma \ref{31709.17h56},
to prove the claim, it suffices to prove:
\beqy\label{3809.21h40} \forall \eps>0,\qquad  \sup_{|\theta|\leq K^{-1}-\eps} \quad\sup_{\substack{1\leq i\leq n\\ n\geq 1}}\sigma_{n,i}<+\infty& \textrm{ and }& \ff{n}\sum_{i=1}^n\sigma_{n,i}^2=1,\\
\label{3809.21h42}  \forall \eps>0,\qquad  \sup_{|\theta|\leq K^{-1}-\eps} \quad\sup_{\substack{1\leq i\leq m_n\\ n\geq 1}}{\sigma'_{n,i}}<+\infty& \textrm{ and }& \ff{m_n}\sum_{i=1}^{m_n}{\sigma'_{n,i}}^2=1.
\eeqy
Note first that   \eqref{3809.21h40} and \eqref{3809.21h42} both hold when $\theta=0$. 
 
A straightforward computation leads, for $\theta\neq 0$, to the formula   \beq \sigma_{n,i}^2
&=&    \f{-16b_n\la_{n,i}^2}{ (b_n^2+\Delta_{n,i}-a_n^2 )^2-4b_n^2{\Delta_{n,i}}} .
\eeq
But (removing the indices)   \beq 
(b^2+\Delta-a^2 )^2-4b^2{\Delta}&=&(2b^2-2ab+4\la^2)^2 -4b^2(b^2-2ab+a^2+4\la^2)\\
&=&16\la^4-16ab\la^2.
\eeq It follows, writing $\gamma_n$ for $\gamma_n(\theta)$, that \be\label{5809.8h29}\sigma_{n,i}^2=\f{b_n}{ a_nb_n-\la_{n,i}^2}=\ff{1+\gamma_n }\times \f{T^{(\f{n}{m_n})}(\gamma_n )}{{T^{(\f{n}{m_n})}(\gamma_n )}-\theta^2\mu_{n,i}^2}=\ff{1+\gamma_n }
\times
\ff{
1-\f{
\theta^2
}{
T^{(\f{n}{m_n})}(\gamma_n  )}
\mu_{n,i}^2
}.\ee
By definition of $\gamma_n(\theta)$, we have $\gamma_n(\theta)\geq 0$, hence ${
T^{(\f{n}{m_n})}(\gamma_n (\theta))}\geq 1$. Since for all $n,i$, $|\mu_{n,i}|\leq K$, it 
follows that 
the first part of \eqref{3809.21h40}
holds. Moreover, by the definition of $\mu_n$ given in \eqref{5809.7h48}, we have $$\ff{n}\sum_{i=1}^n\sigma_{n,i}^2=\ff{1+\gamma_n(\theta)}M_{\mu_n^2}\lf(\f{
\theta^2
}{
T^{(\f{n}{m_n})}(\gamma_n (\theta))}\ri)+\ff{1+\gamma_n(\theta)}.$$
By \eqref{38.09.19h}, it follows that the second  part of \eqref{3809.21h40} also
holds. 

Let us now prove  \eqref{3809.21h42}. A straightforward computation leads, for $\theta\neq 0$ and $i\leq n$, to the formula 
  \beq {\sigma'_{n,i}}^2&=&  \f{a_n}{ a_nb_n-\la_{n,i}^2}.\eeq Hence   for $\theta\neq 0$ and $i\leq n$, \beq {\sigma'_{n,i}}^2&=& \ff{1+\f{n}{m_n}\gamma_n(\theta)}
\times
\ff{
1-\f{
\theta^2
}{
T^{(\f{n}{m_n})}(\gamma_n (\theta))}
\mu_{n,i}^2
},
 \eeq 
  whereas for $i=n+1, \ldots, m_n$, $ {\sigma'_{n,i}}^2=\ff{1+\f{n}{m_n}\gamma_n(\theta)}$. The first part of \eqref{3809.21h42} holds for the same reasons as the first part of \eqref{3809.21h40} above. Moreover,  writing $\gamma_n$ for $\gamma_n(\theta)$, we have $$\ff{m_n}\sum_{i=1}^{m_n}{\sigma'_{n,i}}^2=\f{n}{m_n(1+\f{n}{m_n}\gamma_n )}M_{\mu_n^2}\lf(\f{
\theta^2
}{
T^{(\f{n}{m_n})}(\gamma_n )}\ri)
+\f{n}{m_n(1+\f{n}{m_n}\gamma_n )}+
\f{m_n-n}{m_n(1+\f{n}{m_n}\gamma_n )}
.$$
By \eqref{38.09.19h}, it follows that the second  part of \eqref{3809.21h42} also
holds. 

The proof of  the 
claim  is complete. As a consequence, by \eqref{29.9.10}, the \pro of the event 
$\{P_nD_n^{-1/2}Z_n\in A_n\}$ 
tends to one as $n$ tends to infinity, uniformly on every compact subset of $(-K^{-1}, K^{-1})$ (remember indeed that the matrices $P_n$ and $D_n$ depend on $\theta$).
So by \eqref{5809.8h15},   we have, still writing $\gamma_n$ for $\gamma_n(\theta)$, \be\label{octobre.2010.1} \ff{n}\log(J_n(\theta))=\ee \bes\f{n-m_n}{2n}\log (1+\f{n}{m_n}\gamma_n )-\f{\log (T^{(\f{n}{m_n})}(\gamma_n))}{2}-\ff{2}\int_{t\in [-K,K]}\log \{1-\f{\theta^2}{T^{(\f{n}{m_n})}(\gamma_n)}t^2\}\ud\mu_n(t) +o(1)\ees \bes=-\f{m_n}{2n}\log (1+\f{n}{m_n}\gamma_n )-\f{\log (1+\gamma_n)}{2}-\ff{2}\int_{t\in [-K,K]}\log \{1-\f{\theta^2}{T^{(\f{n}{m_n})}(\gamma_n)}t^2\}\ud\mu_n(t) +o(1).\ees
 By hypothesis, $\mu_n$, which is defined in \eqref{5809.7h48}, converges weakly to $\mu$. 
Using \eqref{5809.7h32} and \cite[Th. C.11]{agz09}, one easily sees that, writing $\gamma$ for $C_\mu^{(\la)}(\theta^2)$, we have  
$$\ff{n}\log (J_n(\theta))= -\f{1}{2\la}\log( 1+\la \gamma) -\f{\log (1+\gamma)}{2}-\ff{2}\int_{t\in [-K,K]}\log \{1-\f{\theta^2}{T^{(\la)}(\gamma)}t^2\}\ud\mu(t) +o(1),$$ where in the case where $\la=0$,  $\f{1}{2\la}\log( 1+\la \gamma)$ has to be understood as $\f{\gamma}{2}$. 
 By \eqref{bird.31709}, one gets, still writing $\gamma$ for $C_\mu^{(\la)}(\theta^2)$, \be\label{14.4.11}I_n(\theta)=\underbrace{\gamma- \f{1}{2\la}\log( 1+\la \gamma) -\f{\log (1+\gamma)}{2}-\ff{2}\int_{t\in [-K,K]}\log \{1-\f{\theta^2}{T^{(\la)}(\gamma)}t^2\}\ud\mu(t) }_{:=f(\theta)}+o(1).\ee 
$I(0)=f(0)=0$ (indeed,  by \eqref{3809.11h39}, $C_\mu^{(\la)}(0)=0$). So to conclude the proof of Theorem \ref{310709.11h}, it suffices to verify that $I$ and $f$ have the same derivatives on $(-K^{-1}, K^{-1})$. Using \eqref{38.09.19h}, it can easily be proved that both derivatives are equal to $\gamma/\theta$.

\noindent {\bf Aknowledgements:} The author would like to thank Myl\`ene Ma\"\i da for her   comments on a draft of this paper and also both referees for their careful reading of the paper and their remarks.

\end{document}